\documentclass[a4paper,12pt]{amsart}

\usepackage{amsthm}
\usepackage{amsmath}
\usepackage{amsfonts}
\usepackage{amssymb}
\usepackage{times}
\usepackage{graphicx}

\theoremstyle{definition}
\newtheorem{defn}{\indent\bf Definition}
\newtheorem{rem}[defn]{\indent\bf Remark}

\theoremstyle{plain}
\newtheorem{lemma}[defn]{\indent\bf Lemma}

\newtheorem{thm}[defn]{\indent\bf Theorem}
\newtheorem{cor}[defn]{\indent\bf Corollary}

\def\tilde{\widetilde}

\def\PSL{\mathop{\rm PSL}\nolimits}
\def\SL{\mathop{\rm SL}\nolimits}
\def\PGL{\mathop{\rm PGL}}

\def\Z{\mathbb Z}
\def\R{\mathbb R}

\def\C{\mathbb C}

\def\L{{\mathcal L}}

\def\H{{\mathcal H}}
\def\A{{\mathcal A}}
\def\S{{\mathcal S}}

\def\D{{\mathcal D}}
\def\B{{\mathcal B}}
\def\K{{\mathcal K}}

\def\cQ{{\mathcal Q}}
\def\gp{{G_p}}
\def\gr{{G_r}}
\def\pg{{\partial \Gamma}}
\def\lc{{\ell^\infty(\Gamma)/c_0(\Gamma)}}
\def\gpt{{\tilde G_p}}

\def\pc{{p_\chi}}
\def\o{{^{\rm op}}}
\def\csl{{C$^\ast$}}
\def\cs{{C^\ast}}

\def\red{{_{\rm red}}}
\def\e{{\epsilon}}

\begin{document}

\title[ Twisted Akemann - Ostrand property]
{ Twisted  Akemann - Ostrand property for ${\rm PGL}_2(\mathbb Z[\frac{1}{p}])$ 
and Ramanujan Petersson Conjectures }
\author[Florin R\u adulescu]{Florin R\u adulescu${}^*$
\\ \\
Dipartimento di Matematica\\ Universit\` a degli Studi di Roma ``Tor Vergata''\\ \\
To Vaughan Jones, in Memoriam}

\maketitle
\baselineskip=20pt

\thispagestyle{empty}

\renewcommand{\thefootnote}{}
\footnotetext{${}^*$ Member of the Institute of  Mathematics ``S.
Stoilow" of the Romanian Academy}
\footnotetext{${}^*$
Supported in part  INdAM-GNAMPA, MIUR Excellence Dept. Project  Dept. Math. Univ  Rome Tor Vergata, CUP E83C180001000., OAAMP,  CUP E81I18000070005, CNCS Romania, PN-III-P1-1.1-TE-2019-
0262, The convex space of sofic representations  }

\begin{abstract}
We prove an extension of  the Akemann - Ostrand 
theorem, regarding the simultaneous,  left and right  regular representations of the free group, modulo compact operators, 
 to the case of the  partial  
action of ${\rm PGL}_2(\mathbb Z[\frac{1}{p}]) \times {\rm PGL}_2(\mathbb Z[\frac{1}{p}])^{\rm op},$  on $\ell^2(\PSL_2(\mathbb Z))$, in the presence of a non-trivial  cocycle.
We use this result and  the operator algebra  techniques developed previously to prove
that the essential spectrum of the Hecke operators is contained in the bounds  
 prescribed by the Ramanujan Peterson Conjectures, which determines the asymptotics of exceptional eigenvalues.
\end{abstract}

\vskip0.5cm

\section*{0. Introduction}

\vspace{0.2cm}

In this paper we prove a variant of the Akemann - Ostrand
theorem (\cite{AO}) 
in the case of the  
action, by partial transformations, by left and and right convolutions, of the group $${\rm PGL}_2(\Z[\frac{1}{p}]) \times {\rm PGL}_2(\Z[\frac{1}{p}])^{\rm op},$$  on the Hilbert space  $\ell^2(\PSL_2(\Z))$.  The  associated unitary representation will be   twisted by  a non-trivial cocycle $\epsilon \in H^2({\rm PGL}_2(\Z[\frac{1}{p}]), \Z_2)$.  We refer to \cite{Klep} for the definition of the twisted group algebra and its rewpresentations.

 In the presence of the cocycle  $\epsilon$  we prove the validity of the Akemann-Ostrand property (we call this twisted Akemann- Ostrand property). It will be   a consequence of the fact that  properly  projective,  unitary representations of   $PSL(2,\mathbb Q_p),$ $ PSL(2,\mathbb R)$ are automatically tempered and hence, as  observed in the paper by Buss, Echeterhoff, Willett (\cite{BE}), the twisted full and reduced C*-algebras associated to these groups coincide. 

The  motivation for the use of the cocycle $\epsilon$ is the fact, described extensively in \cite {ON} (see also \cite{Ra1}, \cite {Ra2}), that the classical Hecke operators acting on Maass forms are  modelled by completely positive maps, that  are obtained by representing, canonically,   the twisted crossed product algebra
\begin{equation}\label{algebra}
 \mathcal A=C^{\ast}(({\rm PGL}_2(\Z[\frac{1}{p}] \times {\rm PGL}_2(\Z[\frac{1}{p}]^{\rm
op}) \ltimes_{\epsilon, \epsilon} C(PSL(2,{\bf Z}_p))),
\end{equation}
into $B(\ell^2(\Gamma))$. In the above formula ${\bf Z}_p$ are the p-adic integers and  $PSL(2, {\bf Z}_p)$ is the  maximal compact subgroup of $PSL(2,\mathbb Q_p)$. 

The cocycle $\epsilon$ appears in a canonical way, because the construction of the completely positive model for the Hecke operators (see Theorem  \ref{Ra}) is based on properties of  the matrix coefficients   of a properly projective, discrete series representation of $PSL(2,\mathbb R)$, restricted to ${\rm PGL}_2(\Z[\frac{1}{p}])$. 
This construction  relies  on the existence of a ciclic trace vector for the restriction of the  representation to the group $PSL(2,\mathbb Z),$ which is  due to the  von Neumann dimension calculations   in \cite{HGJ} and more recently \cite{Jo} .

The above construction (see Paragraph 4) gives a unitary equivalent representation  for the  C$^\ast$- algebra  generated by Hecke operators inside  the image of the canonical representation of the algebra $\mathcal A$ in $B(\ell^2(\Gamma))$. Hence,  to obtain asymptotic estimates of the spectrum for the Hecke operators (i.e estimates of the spectrum modulo the compact operators- that is the essential spectrum) one has to determine the  C$^\ast$-algebra crossed product norm  on the algebra  $\mathcal A$ when represented into $B(\ell^2(\Gamma))$,  modulo  the compact operators. 

The variant of the  Akemann Ostrand property that we prove in this paper shows  that the norm considered above coincides with the  the  C$^\ast$-algebra reduced crossed product norm for the algebra $\mathcal A$. Hence essential norm computations for Hecke operators are transferred in the reduced  C$^\ast$-Hecke algebra, introduced in Bost Connes paper \cite{BC}. Consequently our results prove that  the essential  spectrum estimates are in coincidence with the Ramanujan Petersson estimates, and allow the determination of the asymptotics of exceptional eigenvalues.

The results in this paper have been circulated in the following preprints: arXiv:0802.3548 (\cite{Ra})
 and arXiv:1509.09246. In a recent preprint, by  A. Unterberger,  arXiv:2001.10956 (\cite{unt}), it is proven, with a very different approach,  that  the results in Theorem \ref{thm3} on the essential spectrum of the Hecke operators  are valid, in a stronger form,   for the full spectrum of the Hecke
  operators. 
 
 Acknowledgment. The author is thanking Professor S. Grigorchuk for many remarks on his previous work leading to the proof of the theorems in Section 3. The author is indebted to Professors S. Echterhoff, T. Siebenand, T. de Laat for clarification on the use of the cocycle in Section 2.
 The author is indebted to Professors J. Bassi,  F. Boca,   H. Moscovici,  R. Nest.,  for several discussions regarding topics related to the
subject of this paper. The author is particularly indebted to Professor N. Ozawa  for several comments on this paper and for providing his personal notes (\cite{ON}) for a seminary at the University of Tokyo on the content of the paper (\cite{Ra}). The author is specially thanking to Professor  S. Neshveyev for very pertinent questions on arguments in the proofs.  The author also thanks to Professors Tatiana Smirnova-Nagnibeda and Pierre de la Harpe and the members of the Geometric Group Theory Seminar at the University of G\`eneve where he discussed some aspects related to the constructions in this paper. The author is particularly indebted to Professor Cristophe Pittet for his insight on the use of quasi-rays.  The author thanks Professor Andr\`e Unterberger for his strong  encouragment on this article.

\section{ Overview of The Akemann Ostrand property} 

In this paper we are using the following notational conventions. The symbol $H^{\rm op}$ as an exponent of a group $H$ refers to the fact that we are considering the same group $H$, with the opposite multiplication operation.  If a group $H$ acts by automorphisms or partial automorphisms
on a  locally compact space $X$, we denote by $C^\ast(H\ltimes C_0(X))$ and $C\red^\ast(H\ltimes C_0(X))$ the full and respectively the reduced \csl crossed product algebra (see e.g. \cite{Pe} for definitions). In the presence of  a 2-cocycle 
$\epsilon\in H^2(H,\mathbb T)$, we denote by $C^\ast(H\ltimes_\epsilon C_0(X))$ and $C\red^\ast(H\ltimes_\epsilon C_0(X))$  the corresponding twisted \csl crossed product algebras (\cite{pac}).

Recall that Akemann - Ostrand property (to which we will refer in the sequel as to the AO property) for the free group $F_N$, $N
\geq 2$ asserts (\cite{AO}) the fact  that the $C^{\ast}$ - algebra, generated
in $\B(l^2(F_N))$, by the
$C^{\ast}$ - algebras $C^{\ast}_{\lambda}(F_N)$,
$C^{\ast}_{\rho}(F_N)$ is isomorphic, modulo
the ideal $\K(l^2(F_N))$ of compact operators, to the minimal
$C^{\ast}$ - tensor product $C^{\ast}_{\rm red}(F_N)
\mathop{\otimes}\limits_{\min} C^{\ast}_{\rm red}(F_N^{\rm op}) \cong
C^{\ast}_{\rm red}(F_N \times F_N^{\rm op})$ of the reduced group
$C^{\ast}$ - algebras associated to $F_N$.  The $\cs$-algebras $C^{\ast}_{\lambda}(F_N)$,
$C^{\ast}_{\rho}(F_N)$    are    generated by the  left and respectively, the right
convolution operators with elements in $F_N$ and are isomorphic to $C^{\ast}_{\rm red}(F_N)$. 

The Akemann - Ostrand property has been widely extended by G. Skandalis, who
proved (\cite{Sk}) that the same result remains true for lattices in
semisimple Lie groups of rank 1. Using amenable actions techniques
(\cite{CAD}), Guentner and Higson (\cite{GH}) and later Ozawa (\cite{Oz2}) have further
extended this result, to a  large class of groups, containing the  hyperbolic groups.

The key in Ozawa's approach in proving the AO property for a discrete hyperbolic  group $\Gamma$  is the amenability (see e.g  \cite{CAD}
) of the action of $\Gamma \times \Gamma^{\rm op}$ on the
boundary $\partial(\beta\Gamma)$ of the Stone Cech compactification  of $\Gamma$. This is the Gelfand spectrum of the algebra $\ell^\infty(\Gamma)/c_0(\Gamma)$.
This stronger property for a group $\Gamma$ is called (\cite{Oz2},\cite{BO})  property  
$\S$.

Let
\begin{equation}\label{Koop}
\pi_{\rm
Koop}:C^{\ast}((\Gamma
\times \Gamma^{\rm op}) \ltimes C(\partial(\beta\Gamma)))\rightarrow \B(l^2(\Gamma)).
\end{equation}
be the canonical representation of  the crossed product $C^{\ast}$ - algebra $C^{\ast}((\Gamma
\times \Gamma^{\rm op}) \ltimes C(\partial(\beta\Gamma)))$ into $\B(l^2(\Gamma))$ mapping elements in $\Gamma
\times \Gamma^{\rm op}$ into left and right convolution operators and letting $C(\partial(\beta\Gamma)))\cong \ell^\infty(\Gamma)$ act by multiplication on $l^2(\Gamma)$.

Let 
$$
\pi_{Q} : \B(l^2(\Gamma)) \to \cQ(l^2(\Gamma)) =
\B(l^2(\Gamma))/\K(l^2(\Gamma))
$$
be the projection onto the Calkin algebra (\cite{Ca}).
 The Ozawa's $\S$ property (\cite{Oz}) implies that the representation  
 $$\Pi_Q= \pi_{Q} \circ      \pi_{\rm Koop}:C^{\ast}((\Gamma
\times \Gamma^{\rm op}) \ltimes C(\partial(\beta\Gamma))) \to \cQ(l^2(\Gamma)),$$
\noindent factorises to a representation of the
reduced $C^{\ast}$ - algebra $$C^{\ast}_{\rm red}((\Gamma \times
\Gamma^{\rm op}) \ltimes C(\partial(\beta\Gamma))).$$

\section{ Outline of the results}

We outline the results of this paper. We extend the Akemann - Ostrand property in a specific case, which we describe below. In our setting there is no seemingly analogue for the  property $\S$. 
Recall  that by the results in \cite{Sa}, the group $G = \PGL_2(\Z[\frac{1}{p}])$, $p$ a prime $\geq 2$, does not have the property   $\S$.


  We denote, respectively, by $G_p, G_r$, $\Gamma$ the  groups $PSL(2,\mathbb Q_p)$,  $ PSL(2,\mathbb R)$ and $PSL(2,\mathbb Z)$. Let $K= PSL(2,{\bf Z}_p) $ be the maximal compact subgroup in $PSL(2,\mathbb Q_p)$. Here ${\bf Z}_p$ are the $p-$adic  integers.  Let 
   $\tilde G_p$, $\tilde G_r, \tilde K$, $\tilde \Gamma$ be the respective  $SL_2$ counterparts of the previous groups. Note that second series of groups are a central extensions, by $\mathbb Z_2$, of the previous ones. 
  
  Let $\chi$ be the non-trivial character of $\mathbb Z_2$, Let $s$ be the diagonal matrix with entries -1 on the diagonal, which is a central element of each of those groups and let
  \begin{equation}\label{pchi}
   p_\chi=
 \frac{1-s}{2}.
 \end{equation}
Then $p_\chi$ is a selfadjoint   projection   which belongs to the C$^*$-algebra (or multiplier algebra) of each of the previous groups. Let $\epsilon$ be the (non canonical) 2 cocycle with values in $\Z_2$, corresponding to the projective representation $\pi_{13}$ in the analytic series of $\gr=PSL(2,\mathbb R)$ (see e.g. \cite {HGJ}).
 
We consider the twisted C$^*$ groupoid  crossed product algebra (we refer to \cite {JR} for definitions)
 $$\mathcal A=C^{\ast}((G \times G^{\rm
op}) \ltimes_{\epsilon, \epsilon} C(K)).$$

The above crossed product is  a grupoid crossed product, since  the action 
of $G \times G^{\rm
op}$ by left and right multiplication, on $K$, is a partial action (see the next  section).
It is easy to see that the  algebra $\mathcal A$ is  a corner of a (genuine) twisted C$^\ast$- crossed product algebra.  In fact, let

$$\mathcal E= C^{\ast}((G \times G^{\rm
op}) \ltimes_{\epsilon, \epsilon} C_0(G_p)).$$
 Let $\chi_K\in C_0(G_p)  $ be the charcteristic function of $K$.  Then it is obvious to check that
\begin{equation}\label{corner}
\mathcal A\cong \chi_K \mathcal E\chi_K.
\end{equation}

The algebra $\mathcal A$ has a canonical representation $\pi_{\rm Koop}$ into $B(l^2(\Gamma))$ which we describe below. To avoid the complication 
of having to use   the  (non-canonical) 2-cocycle $\epsilon$, and the twisted crossed product, we use the technique in \cite{pac} and identify $l^2(\Gamma)$ with 
$p_\chi l^2(\tilde \Gamma)$. Let
$$\mathcal {\tilde A}=C^{\ast}((\tilde G \times \tilde G^{\rm
op}) \ltimes C(\tilde K)).$$
Then, as in formula  (\ref{Koop}),  there is a canonical representation 
\begin{equation}\label{kooptilde}
\tilde{\pi_{\rm Koop}}: \mathcal {\tilde A}\rightarrow B(\ell^2(\tilde \Gamma)).
\end{equation}
In the above formula, the action of $C(\tilde K)$ on $l^2(\tilde \Gamma)$  is determined by letting the characteristic function of a coset of a subgroup $\tilde K\cap g\tilde Kg^{-1}$ act on $l^2(\tilde \Gamma)$  by multiplication with the characteristic function of $\tilde \Gamma\cap g\tilde \Gamma g^{-1}$, $g\in G$.
Then we define
\begin{equation}\label{inc}
\pi_{\rm Koop}=p_\chi\tilde{\pi_{\rm Koop}}:  \mathcal A=\pc \mathcal {\tilde A}\pc\rightarrow B(p_\chi \ell^2(\tilde \Gamma))\cong B(\ell^2(\Gamma)).
 \end{equation} 
 

    The main theorem of the second section is
 \begin{thm}\label{t1}
 Let $\Pi_Q$ be the representation of the algebra $\A$ into $Q((l^2( \Gamma)))$ obtained by composing the representation $\pi_{\rm Koop}$ of $\mathcal A$ into $B(\ell^2(\Gamma))$  from formula (\ref{inc}) with the projection $\pi_Q$ onto the Calkin algebra. 
 Then $\Pi_Q$ factorises to the reduced crossed product algebra $$\mathcal A_{\rm red}=C^{\ast}_{\rm red}((G \times G^{\rm
op}) \ltimes_{\epsilon, \epsilon} C(K)).$$
 \end{thm}
 To prove the theorem we will first prove
that
 
  \begin{thm}\label{t2}
 The (maximal) crossed product algebra 
 $$\mathcal B=C^{\ast}(( G \times  G^{\rm
op}) \ltimes_{\epsilon, \epsilon} C_0(G_r\times K \times \gr\o ))$$  coincides with the reduced
crossed product algebra
$$\mathcal B_{\rm red}= C^{\ast}_{\rm red}((G \times   G^{\rm
op}) \ltimes_{\epsilon, \epsilon} C_0(G_r\times K \times \gr\o )).$$
Moreover the same holds true if $G_r\times K \times \gr\o$ is replaced by
$P^1(\mathbb R)\times K\times P^1(\mathbb R), $ where $P^1(\mathbb R)$ is the projective plane.
\end{thm}
Theorem \ref{t2} will be used to prove Theorem \ref{t1} by showing that inside $\ell^\infty(\Gamma)/c_0{\Gamma}$ there is a $G \times  G^{\rm
op}$-equivariant copy of the abelian $C^\ast$-algebra $C_0(G_r\times K \times \gr\o )$.

We then use Theorem \ref{t1} to find estimates on the Hecke operators.  We proved in \cite{Ra} (see also \cite{ON}) that the classical Hecke operators (\cite{He})  $T_{p^n}$, $n\geq 1$,  acting on Maass forms on the upper halfplane $\mathbb H$, associated to 
$\sigma_{p^n}= \left(\begin{array}{cc} p^n & 0 \\ 0 & 1 \end{array} \right)$, $n$ a positive integer, are unitarely equivalent
(up to a commuting phase that is a function of the invariant laplacian) to the image  in the representation of the algebra $\mathcal A$  defined  in formula (\ref{inc}) of an element of the form
$$S([\Gamma\sigma_{p^n}\Gamma])=\chi_K(t([\Gamma\sigma_{p^n}\Gamma])\otimes 
t([\Gamma\sigma_{p^n}\Gamma]))\chi_K\in\mathcal A\cong \chi_K \mathcal E\chi_K.$$
In the previous formula we use the identification from formula (\ref{corner}). Moreover $$t([\Gamma\sigma_{p^n}\Gamma]\in \cs\red(G,\epsilon))$$ is the image of the coset $[\Gamma\sigma_{p^n}\Gamma]$, viewed as  an element of the Hecke algebra $\mathcal H_0$ of double cosets of $\Gamma$ in $G$, via a canonical representation $$t:\mathcal H_0:\rightarrow \cs\red(G,\epsilon).$$
This representation  is constructed in \cite {Ra} using the matrix coefficients of the representation $\pi_{13}$ (see Section 4 for more details).

Using the above mentioned  results from (\cite{ON}, see also \cite {Ra}) we prove

\begin{thm}\label{thm3}
The essential  spectrum of the Hecke operator  $T_p$,   acting on Maass forms on the upper halfplane $\mathbb H$ (i.e. the spectrum after modding out compact operators) is contained in the interval $[-2\sqrt p, 2\sqrt p]$. 

Consequently, the possible exceptional eigenvalues of the Hecke operator $T_p$  necessarily accumulate to  the endpoints of the   interval  $[-2\sqrt p, 2\sqrt p]$, predicted by the Ramanujan Petersson conjectures.
 \end{thm}

\section{The partial action of  ${\rm PGL}_2(\Z[\frac{1}{p}]) \times {\rm PGL}_2(\Z[\frac{1}{p}])^{\rm op},$} 

It is well
known (\cite{He}), that $\Gamma$ is almost normal in $G$. The almost normal property for the subgroup $\Gamma$ of $G$ signifies  that  for all
$g \in G$ the subgroup 
\begin{equation}
\Gamma_g = g\Gamma g^{-1} \cap \Gamma \subseteq \Gamma,
\end{equation}
 has
finite index $[\Gamma : \Gamma_g]$.

The group $G$ acts naturally, by conjugation, by partial isomorphisms, on $\Gamma$.
Indeed for $g \in G$, the conjugation by $g$ on $G$, will restrict to
a partial isomorphism
\begin{equation}\label{delta}
\Delta(g) : \Gamma_{g^{-1}} \to \Gamma_g.
\end{equation}

It is well known (see e.g. \cite{Serre}) that for   the modular group $\Gamma= PSL(2,\mathbb Z)$, we have that $[\Gamma : \Gamma_g] = [\Gamma
: \Gamma_{g^{-1}}]$, for all $g \in G$. We consider the family of
maximal normal subgroups $\Gamma_g^0$ contained in $\Gamma_g$.

Let $K$ be the compact space obtained as the inverse limit of the
finite coset spaces $\Gamma / \Gamma_g^0$ as $g \to \infty$. Then $K$
is the totally disconnected subgroup $PSL(2,{\bf Z}_p)$. The   Haar measure $\mu_K$ on $K$ defined by the
requirement that the compact set corresponding to the closure of a  coset $s \Gamma_g, s\in\Gamma$
in the profinite topology, has Haar measure equal to 
$\displaystyle\frac{1}{[\Gamma : \Gamma_g]}$, $g \in G$.

The condition that $[\Gamma : \Gamma_g] = [\Gamma : \Gamma_{g^{-1}}]$,
implies that the partial transformation $\Delta(g)$, introduced in formula (\ref{delta}), induced by
conjugation with $g \in G$, preserves the Haar measure $\mu_K$ on $K$.

There is a natural action of $G \times G^{\rm op}$ on $K$.
An element $(g_1, g_2) \in G \times G^{\rm op}$ acts by partial
transformations on $K$, by mapping,  $k \in K$ into
$g_1 kg_2^{-1}$, if the later element also belongs to $K$. Thus, the domain of $(g_1, g_2
)$, as a partial
transformation on $K$, is
$$
\D_{(g_1, g_2)} = \{ k \in K \mid g_1Kg_2^{-1} \in K \} = K \cap
g_1^{-1}Kg_2 = K \cap g_1^{-1}Kg_1(g_1^{-1}g_2).
$$

We use the notation $K_g = K \cap gKg^{-1}$. In our construction this   is the
closure, in the profinite completion, of $\Gamma_g$. 




\section{Proof of Theorem \ref{t2}}

To prove theorem 2 we will consider first the  action by left and right actions of the groups $\gp, \gp^{\rm op}$ on $C_0(G_p)$. Consider the algebra 
\begin{equation}\label{ac}
\tilde{\mathcal C}=[p_\chi \otimes (p_\chi)\o]  [C^{\ast}((\gpt \times   \gpt^{\rm
op}) \ltimes C_0(\gpt))].
\end{equation}

Here $p_\chi$ is the projection defined in formula (\ref{pchi}) and   $(p_\chi)\o$ is the projection corresponding to the generator $s$ of $\mathbb Z_2$ viewed as an element in
$\gp^{\rm 
op}$ Clearly then  $p_\chi \otimes (p_\chi)\o$ is a central projection in the multiplier algebra.

In the statement of formula \ref {inc} we used instead the projection $\pc$. Since $\pc$ is a central projection in the corresponding group algebra, it will be no difference in the representation from  formula (\ref{inc}) if we use any of the two projections ($\pc$ or $p_\chi \otimes (p_\chi)\o$),  as they coincide in the given representation.

We choose a Borel lifting  of $\gp$ to $\gpt$. This choice induces a (non-necessary canonical) two cocycle $\tilde \epsilon$ with values in $\mathbb Z_2$  on $G_p$. We may assume, with no loss of generality, that the restrictions $\tilde \epsilon|_G$ and $\epsilon|_G$ coincide.  Then
the algebra $\tilde{\mathcal C}$ is isomorphic to the twisted crossed product algebra 
$$\mathcal C=C^{\ast}((\gp \times   \gp^{\rm
op}) \ltimes_{\tilde\epsilon, \tilde\epsilon}C_0(\gp)).$$

We prove first the following lemma

\begin{lemma}\label {l1}
The full crossed product C$^*$- algebra $\tilde{\mathcal C}$ defined in formula (\ref {ac}) coincides with the reduced crossed product algebra
$$\tilde{\mathcal C}_{\rm red}=[p_\chi \otimes (p_\chi)\o]  [C^{\ast}_{\rm red}((\gpt \times   \gpt^{\rm
op}) \ltimes C_0(\gpt))]\cong C^{\ast}_{\rm red}((\gp \times   \gp^{\rm
op}) \ltimes_{\tilde\e, \tilde\e} C_0(\gp)) .$$
\end{lemma}
\begin{proof}
To simplify the notations in the proof,  we use   the $\epsilon$ terminology.  
Thus we have 
\begin{equation}\label{dc}
\tilde{\mathcal C}\cong  C^{\ast}((\gp \times   \gp^{\rm
op}) \ltimes_{\tilde\e, \tilde\e} C_0(\gp)).
\end{equation}
Let $\cs(\gp,\tilde\e)$ be the twisted group \csl-algebra of the group $\gp$ with respect the cocycle $\tilde\e$.
This is isomorphic to $\pc\cs(\gpt)$. We use the approach from the paper \cite{BE}, Proposition 5.24. Because in the supplementary series of irreducible unitary representations of $\gpt$, the element $s$ is always mapped into the identity (see chapter 2.3.7 of \cite{GG} and also \cite{Sal}, Section 14) it follows that
$$\cs(\gp,\tilde\e)=\cs\red(\gp,\tilde\e).$$

The crossed product algebra $\tilde{\mathcal C}$ may be written as the iterated crossed product 
\begin{equation}\label{dcc}
\cs(\gp^{\rm
op}\ltimes_{\tilde\e} [\cs(\gp\ltimes_{\tilde\e} C_0(\gp))]).
\end{equation}
By the amenability of the action of $\gp$ on $\gp$ it follows that the inner crossed product algebra 
$\cs(\gp\ltimes_{\tilde\e} C_0(\gp))$ coincides with the reduced crossed product algebra 
$\cs\red(\gp\ltimes_{\tilde\e} C_0(\gp))$ and this in turn coincides with the C$^\ast$-algebra
$\mathcal K(L^2(\gp, \nu_p))$, where $\nu_p$ is the Haar measure on $\gp$. Because of the observation at the start of the proof and because of the Proposition 5.24  in \cite{BE} it follows that also the outer maximal \csl-crossed product in formula (\ref{dcc}) coincides with the reduced \csl crossed product.
Since both iterated \csl-crossed products in formula (\ref{dcc}) coincide with the reduced \csl-crossed products, it follows that the algebra $\tilde{\mathcal C}$ in formula (\ref{dc}) coincides with the reduced crossed product algebra
$$C^{\ast}\red((\gp \times   \gp^{\rm
op}) \ltimes_{\tilde\e, \tilde\e} C_0(\gp)).$$

\end{proof}

The next lemma combines the left right action of $\gp$ and $\gr$.

\begin{lemma}\label{g} Let $\mathcal G=(\gr\times\gr\o)\times (\gp\times\gp\o)$ act on 
\begin{equation}\label{x}
\mathcal X=\gr\times(\gr)\o\times \gp. 
\end{equation}
Here $\gp\times\gp\o$ acts by left and right action on
$\gp$, and $\gr$ and  respectively $(\gr)\o$ act on $C_0(\gr)$ and respectively $C_0((\gr)\o)$.
Then the twisted full \csl-crossed product 
$$\cs(\mathcal G\ltimes_{\tilde\e,\tilde\e} C_0(\mathcal X))$$ coincides with the twisted reduced crossed product
$$\cs\red(\mathcal G\ltimes_{\tilde\e,\tilde\e} C_0(\mathcal X)).$$

  \end{lemma}
  \begin{proof} The crossed product algebra $\cs(\mathcal G\ltimes_{\tilde\e,\tilde\e} C_0(\mathcal X))$ splits as
  \begin{equation}\label{tens}
  \cs(\gr\ltimes C_0(\gr))\otimes   \cs(\gr^\o\ltimes C_0(\gr^\o))\otimes
   \cs((\gp\times\gp\o)\ltimes_{\tilde e, \tilde e}C_0(\gp)).
   \end{equation}
   
  By the amenability of the action of $\gr$ on $\gr$, the first two algebras are nuclear and hence  all the tensor products are minimal tensor products. The third factor in the tensor product, by the previous lemma, coincides with the reduced tensor product.

  \end {proof}
  
  \begin{cor}\label{p1} Let $P$ be the quotient $\gr/AN\cong P^1(\mathbb R)$, endowed with the standard action of $\gr$.
  Let $\mathcal Y=P\otimes \gp\otimes P^\o$. Let $\mathcal G$ act on $\mathcal Y$ as in the previous statement.

  Then the twisted \csl-crossed product 
$$\cs(\mathcal G\ltimes_{\tilde\e,\tilde\e} C_0(\mathcal Y))$$ coincides with the reduced crossed product
$$\cs\red(\mathcal G\ltimes_{\tilde\e,\tilde\e} C_0(\mathcal Y)).$$

  \end{cor}
  \begin{proof} The proof of the previous statement was using only the amenability of the action of $\gr$ on $\gr$. Hence the same proof works if replace $\gr$ by $P$.
  
  \end{proof}

  To prove Theorem  \ref{t2}  we use the following lemma. We are restricting he to the  use of the group $\mathbb Z_2 \oplus \mathbb Z_2$, but it is very likely this works in a more general setting.
  
  \begin{lemma}\label{downlattice}
  Let $\tilde H$ be a locally compact, separable,  exact group. Assume that the center of 
  $\tilde H$ is equal to a copy of $\mathbb Z_2 \oplus \mathbb Z_2$, with generators $s_1,s_2$. Let $p_i=\frac {1-s_i}{2}, i=1,2$ be the projections corresponding to the eigenvalue 1 in the center of the group algebra of $\tilde H$   and let $p=p_1p_2$.
  Let $X$ be a locally compact separable space, endowed with a continuous action $\alpha$ of the group $\tilde H$ and assume that the central elements $s_1,s_2$ act trivially on $X$.
  Assume that 
  $$p\cs_{\rm max}(\tilde H\ltimes C_0(X))=p\cs_{\rm red}(\tilde H\ltimes C_0(X)).$$
  Let $\Delta$ be a countable discrete subgroup of $\tilde H$, containg $s_1,s_2$. Then
  $$p\cs_{\rm max}(\Delta\ltimes C_0(X))=p\cs_{\rm red}(\Delta\ltimes C_0(X)).$$
  
  \end{lemma}
  
     \begin{proof}We will use an adapted version of the proofs in Chapter 5 of 
     \cite {BE}. This is  a twisted version of the fact that in the commutative case, the weak containment property (\cite{BE}) implies strong amenability (topological amenability) in the sense of \cite{CAD}. Let $A=C_0(X)$ and  let $A''_{\alpha, \epsilon}$ be the image of the algebra $A$ in the "restricted" universal covariant  representation $\pi_{u,\epsilon}$ of the system $(C_0(X)),\alpha)$. The word "restricted" has the meaning that we only consider covariant representations of the system that verify the condition $\pi(s_i)= -1, i=1,2$. The subscript $\epsilon$ is used here just for notational purposes to denote the "restricted" case. (If we replace the group $\tilde H$ by the quotient $H=\tilde H\backslash( \mathbb Z_2 \oplus \mathbb Z_2)$ then we get a twisted crossed product depending on a two cocycle on $H$ with values in  $\mathbb Z_2.$)
     
     Clearly if $A''_\alpha$ is the enveloping $\tilde H$-von Neumann algebra (\cite{Ik}, see also the definition in \cite{BE}), then the projection $p$ commutes with $A''_\alpha$ and
    $A''_{\alpha, \epsilon}= p  A''_\alpha$. The proof of  commutant amenability (see \cite {BE}) for definitions )  for the "restricted" universal covariant representation follows exactly the lines in the above mentioned paper.
    
     Indeed any covariant representation $\pi $ of the system $(A,\alpha)$ into $B(H)$, for some Hilbert space $H$ that verifies the condition $\pi(s_i)=-1. i=1,2$ is continuous, by the  hypothesis, with respect to the reduced crossed product and hence, by the injectivity of the reduced crossed product  and by Arveson Theorem,  $\pi$ extends to a completely positive map defined on
    $ p\cs_{\rm red}(\tilde H\ltimes (C_0(X) \times C_{\rm ub}(\tilde H)))$.
    Here  (as in \cite{OS}, \cite {BK}, \cite{BE})  $C_{\rm ub}(\tilde H)$ is the algebra of bounded, left uniformly continuous with respect $\tilde H$,  $\mathbb C$-valued,  continuous functions on $\tilde H$. 
    
    Arguing as in Lemma 5.6 in \cite{BE}, by looking at the multiplicative domain of the unital completely positive map, we obtain a $\tilde H$ equivariant unital completely positive map from $ C_{\rm ub} (\tilde H)$ into $\pi(A)'$. Using the amenability of the action of $\tilde H$ on $ C_{\rm ub} (\tilde H)$ (proved in Proposition 2.5 in \cite{OS}), since  $\tilde H$ is exact,   we obtain as in \cite {BE} a net of $\pi(A)'$ valued positive definite functions defined on $\tilde H$, converging to 1 ultraweakly and uniformly on compact subsets of $\tilde H$.
    
   Letting $\pi$ be the "restricted" universal representation  $\pi_{u, \epsilon}$ and using as in \cite {Mats} (see also Section 5.2 in \cite {BE}) the Haagerup standard form (\cite{Ha}) for the von Neumann algebra $A''_{\alpha, \epsilon}$ we get the  (von Neumann) amenability  of the action of $\tilde H$ on $Z(A''_{\alpha, \epsilon})=pA''_{\alpha}$, that is we get  a net $(\theta_i:\tilde H\rightarrow pZ(A''_{\alpha})=pA''_{\alpha})_{i\in I}$, of norm continuous, compactly supported functions of positive type such that $||\theta_i(e)||\leq 1, i\in I$  and $\theta_i(h)\rightarrow  p$, for $h\in \tilde H$, ultraweakly and uniformely on compact subsets of $\tilde H$ (see Definition 1.3 in \cite {BE}).  
   
   The arguments in
  Corollary 4.14 in  (\cite {BeCr}) prove that one can improve the  statement above so as to assume that there exists a choice for the  functions $\theta_i$ such that they take values into $pA$.
   The fact that we can extend the results in \cite{BeCr} to the present (twisted case-  i.e the central projection $p$ substituting 1)   is because  approximation arguments in the proof of Proposition  4.6 in the above mentioned paper,  are based on the variant of the Kaplansky density Theorem in Corollary 2.7 in \cite{Ze} which is still valid in this case.
   
  Restricting the positive definite functions $\theta_i$ to $\Delta$, and using the fact that $\Delta$ is a discrete subgroup of $\tilde H$, we obtain the statement of the Lemma.
    
    \end{proof}

  With this we can conclude the proof of Theorem \ref{t2}.
  \begin{proof}{(Theorem \ref {t2})}
  We use the fact the group  $G$ is a lattice in $G_{pr}=\gr \times \gp$. 
  
  We apply the previous lemma with $\tilde H= H_1\times \tilde H_2$ and  $H_1=\gr\times\gr^{\rm op}$, $\tilde H_2=\tilde\gp\times \tilde \gp^{\rm op}$, and $X=X_1\times X_2$ with 
  $X_1=\gr\times(\gr)\o,$  (or $P\times P\o)$ and $X_2= G_p$, $\Delta={\rm Sl}(2, \mathbb Z[1/p])\times {\rm Sl}(2, \mathbb Z[1/p])\o$ and use Lemma \ref{g} and Corollary \ref{p1}.  The algebra $C(X)$ is being identified to the continuous functions on $H_1\times \tilde H_2$ that in the second variable are invariant to the action of $(s_1,s_2)$.   The group $\tilde H$ is exact by the results in \cite{GHW}.  As explained in the introduction multiplying the crossed product  C$^\ast$-algebras with the projection $[p_\chi \otimes (p_\chi)\o]$ yields the twisted crossed product  C$^\ast$-algebras.  
  
   To obtain the analogue result also  for the crossed product algebra obtained by   having $C(K)$   replace $C_0(G_p)$, it is sufficient to observe that the property of having the full C$^\ast$ algebra crossed product equal to te reduced crossed product passes to corners (the algebra obtained multiplication by a projection to the left and to the right of the algebra). It  is thus sufficient to use the observation contained in formula (\ref{corner}).

  \end{proof}

  \section{ a $G\times G\o$ equivariant factor of $\partial(\beta\Gamma)$ that is isomorphic to $\pg\times K\times\pg\o$ }
  
  Let $\partial \Gamma$ be the (Gromov) boundary of $\Gamma$ (which is independent on the choice of generators, see e.g.  \cite {kks}). Let $S,T$ be the standard system of generators of $\Gamma\cong \Z_2\ast \Z_3$ and let
  $\nu$ be the homogenous  measure on $\partial \Gamma$ associated  to this system of generators. 
  Let $\pi_d$ be the $\Gamma$ equivariant projection from $\partial \Gamma$ onto $P^1(\mathbb R)$ constructed in \cite {spiel}, \cite{series} (the Spielberg disconnection). In this section we describe first  a canonical $G$ action on $\partial \Gamma$ that is 
  $G$-equivariant. We construct an embedding of   $C_0(\partial \Gamma\times \partial \Gamma\o)$ inside $\ell^\infty(\Gamma)/c_0(\Gamma)$, and prove that the corresponding embedding is $G\times G\o$ equivariant with respect to the canonical action of $G\times G\o$ on
  $\ell^\infty(\Gamma)/c_0(\Gamma)$.
  
The construction introduced  in this chapter has been recently realised in a different way, for groups like $SL(3,\mathbb Z)$, in \cite{BR}, based on the theory of quasi-projective transformations defined by Furstenberg (\cite {F1}, \cite {F2}). In the case of $\Gamma=PSL(2,\mathbb Z)$  this is easy described: one observes that there is a canonical  $G\times G^{\rm op}$-equivariant map from $\partial(\beta\Gamma)$ onto the compact space of 2 by 2 matrices  of rank 1 (modulo the scalars). The equivariant map onto $P^1(\mathbb R) \times  P^1(\mathbb R)$ is then obtained by composing with the map that to a linear transformation associates the range and corange of the transformation. 
  
  We divide the proof of the construction   in several steps.
  
  \begin{lemma}\label{action}
  There exists a canonical extension of the action of $\Gamma$ on $\pg$ to an  action of $G$ on $\pg$. This action  is $G$-equivariant with respect to $\pi_d$.
  \end{lemma}
  
  \begin{proof} Every $g$ in $G$ defines a partial action on by conjugation with domain $\Gamma_{g^{-1}}=g^{-1}\Gamma g\cap\Gamma$ and  codomain $\Gamma_{g}=g\Gamma g^{-1}\cap\Gamma$. Since both subgroups $\Gamma_{g^{-1}}$ and $\Gamma_{g}$ have finite index in 
  $\Gamma$, their Gromov boundaries coincide with $\pg$ (see \cite {gkn}).
  Moreover as the boundary is independent of the chosen system of generators (see e.g. \cite {kks}) and hence  the proof in \cite {cooper} extends to this situation and it implies   that the  conjugation action by $g\in G$  gives a continuous automorphism of $\pg$ that clearly extends the action of $\Gamma$.  
  Recall that  by $\lambda_g$ respectively $\rho_g$ we denote  the left (respectively right multiplication) by $g$ on $l^2(\Gamma)$.
    
  
  Recall that in \cite{spiel}, page 781 (see also \cite{series}) the projection $\pi_d$ is constructed as follows.
  One fixes a point $z\in \mathbb H$ and if $\zeta=\alpha_1\alpha_2...$ is the infinite word expression in terms of the generators then 
  \begin{equation}\label{limit}
  \pi_d(\zeta)=\lim_{n\rightarrow \infty}\alpha_1\alpha_2...\alpha_n (z),
  \end{equation}
  The above definition is independent of $z$ (\cite{series}).

  An alternative  method to describe the action by conjugation of  $g\in G$ on the boundary $\partial \Gamma$ is as follows: Consider a sequence of words $w_1,w_2, ..., w_n,...$ converging to $\zeta\in\pg $. Then there is a constant $M$ such that for any $n$ we find an element $w'_n$ in $\Gamma_g$ that is at distance less than $M$ from $w_n$. We than consider the sequence $w'_n=gw_ng^{-1}$.

  By \cite{raney} the conjugacy action of $G$ on words in $\Gamma$ is that of a Turing machine. Hence the  conjugation,  by the arguments in \cite {cooper} eventually, when $n \rightarrow \infty$, preserves the initial part     of the  words $w'_n$. Since every automorphism of a finitely generated group  is Lipschitz with respect the group metric, it follows that  the sequence $w'_n$  converges to some element $\zeta'\in\pg$   belonging to the boundary. Then $\zeta'$ is the result of the conjugation action by $g$ on $\zeta$.  Here we are using the fact  that the sequences  $w_n, n\in\mathbb N$, and
   $gw_ng^{-1}, n\in\mathbb N$ are defining quasi-rays, and for a hyperbolic space, the Gromov boundary is  defined (\cite{GH}, Chapter 7, Proposition 4) as the equivalence class of quasi-rays (modulo the equivalence relation defined by being at finite distance).

 Let  $z'$ be equal to $\pi_d(\zeta')$.
 Then 
 $$\pi_d(\zeta')=\lim_{n\rightarrow \infty}w'_n (z)=g\lim_{n\rightarrow \infty}w_n (g^{-1}z).$$
    Since  the formula \ref{limit}  doesn't depend on the choice of $z$, it follows that $z'=gz$ and hence it follows that $\pi_d$ is $G$-equivariant.
  
    \end{proof}
  
 \begin{rem}
 An intuitive  method   to define the action of $g\in G$ on $\pg$ is to consider
 along $\zeta\in \pg$ a tube $T$ of fixed width (greater  than the length of coset representatives  for $\Gamma_{g^{-1}}$.  Then $T\cap \Gamma_{g^{-1}}$ is infinite, and $gTg^{-1}$ will be contained in an infinite tube, of eventually larger width, that unequivocally defines the action of $g$ on $\zeta$. The same happens (i.e. $gTg^{-1}_1$ is contained in a larger tube defining the same element on the boundary) if instead we use an element $g_1\in G$ such that $g\Gamma g^{-1}_1\cap \Gamma$ is non-void
 \end{rem}

  \begin{lemma}
 The canonical embedding of 
  $C_0(\pg)$ inside $\lc$ has the property that the commutative C$^\ast$ algebra $A$ generated by $C_0(\pg)$ and $C(K)$ inside $\lc$ is isomorphic to $C_0(\pg\times K)$. Here $C(K)$ is realised as the algebra generated by characteristic functions of cosets of subgroups of the form $K_g$, $g\in G$.
  
  Moreover the partial  action of $G$ by conjugation on $\lc$ leaves the algebra  $A$ invariant and on the factor $C_0(\pg)$ it induces the  action constructed in Lemma \ref{action}, while on $C(K)$ it induces the partial action by conjugation, as in formula (\ref{delta}).
  
  \end{lemma}
  
  \begin{proof} 
  
  One considers the commutative \csl- algebra B in $\lc$ generated by characteristic functions of sets of the form:
  \begin{equation}\label{starts}
  A_w=\{w_1\in \Gamma  \ | \  w_1 {\rm \ starts \ with\  } w\}, \ w \in \Gamma.
  \end{equation}
  
  Then $A$ is the commutative \csl-algebra generated by $B$ and $C(K)$ in $\lc$.
  To prove that the algebra $A$ is a faithful realisation of the commutative 
  \csl algebra $C(\pg\times K)$ , we need to exhibit a state (measure) on $\lc$ whose pushback to $A$ is a product state. To do this let $\omega$ be a free ultrafilter on $\mathbb N$. For $n\in \mathbb N$, let $B_n$ be the ball of radius $n$, with respect to the generators $S,T$ in $\Gamma$
  and let $\mu_\omega$ be the Loeb (\cite {Lo}, \cite{Li}) counting measure associated to the sets $(B_n)_{n\in \mathbb N}$ and the ultrafilter $\omega$.
Consider the pushback from the Loeb measure space associated to the counting measure to the Stone-Cech compactification $\beta \Gamma$ and consider the further pushback of the Loeb measure to $\pg\times K$. It is tautological that the marginal of the pushback measure on  $\pg$ is the homogenous measure $\nu$ (in the terminology introduced in the paper \cite{kks}, with respect to the generators use the construct the balls $B_n, n\in \mathbb N$).
  
  The results in \cite{Gr} prove that for a coset $C\subseteq K$ of a group of the form $K_g$, $g\in G$,  we have that
  $$\lim_{n\rightarrow\infty} \frac {{\rm \ card\ }(B_n\cap C) }{{\rm \ card\ }(B_n)}=\frac{1}{[\Gamma:\Gamma_g]}.$$
  Hence the marginal of the pushback of the Loeb measure on $K$ is the invariant Haar measure.  Thus the algebra $A$ is the algebra $C(\pg\times K)$.
  
  For $g\in G$, the continuity of the homeomorphism $g$ on $\pg$ shows that $g$ will map (modulo finite subsets) sets of  the form $A_w\cap \Gamma_g$ into a reunion of such sets, and hence the conjugation action by $g$ on $\lc$ invariates the algebra $A$. By the previous Remark, the action has the form described in the statement on the two factors.

  \end{proof}
  
  In analogy with formula (\ref{starts}), for $w\in\Gamma$, let 
  
  \begin{equation}\label{ends}
  A^w=\{w_1\in \Gamma  \ | \  w_1 {\rm \ ends \ with\  } w\},.
  \end{equation}

  \begin{thm}\label{extension}
   Let $G\times G\o$ act on $\pg \times K\times \pg\o$ by letting $G$ act on $\pg$ as in Lemma \ref{action}, and trivially on $\pg\o$. Similarly let $G\o$ act on $\pg\o$ as in Lemma \ref{action}, and trivially on $\pg$. 
  We let $G\times G\o$ act by partial isomorphisms on $K$.
  Let $C$ be the commutative \csl subalgebra of $\lc$ generated by
  characteristic functions of the form $A_{w_1}$, $A^{w_2}$, $w_1,w_2\in\Gamma$ and characteristic functions of cosets in $\Gamma$ of subgroups of the form $\Gamma_g$, $g\in G$.
  
   Then $C$ is isomorphic to  $C(\pg\times K\times\pg\o)$ and the partial action of $G\times G\o$ on $\lc$ invariates $C$ and induces the action described above.

  \end{thm} 
  
    \begin{proof}
   We prove first  that the algebra $C$ is isomorphic to $C(\pg\times K\times\pg\o)$. Let $w_1, w_2$be two   words in $\Gamma$ and $g_1, g_2 $ in $G$.   Then, modulo finite sets, we have that the set
   $$A_{w_1}\cap A^{w_2}\cap g_1\Gamma g_2^{-1}$$
   is the finite union, modulo a finite set of elements,  of the sets of the  form, 
   $$ A_{w'_1}\cap A^{w'_2}\cap g_1\Gamma g_2^{-1}$$
  where $w'_1$, respectively $w'_2 $ are successors, at distance 1, of $w_1$ and respectively $ w_2$. The  above intersections  are always infinite, and since
   the following   equalities, valid in $\ell^\infty(\Gamma)/c_0(\Gamma)$,  (with $w_1,w_2, w'_1,w'_2, g_1,g_2$ as above)
   $$\chi_{A_{w_1}\cap A^{w_2}\cap g_1\Gamma g_2^{-1}}=\sum \chi_{A_{w'_1}\cap A^{w'_2}\cap g_1\Gamma g_2^{-1}},$$
    describe  completely the algebra $C$, the first part of the statement is thus proved.

 The action induced on the algebra $C$, in the Calkin algebra representation, by an element $(g_1,g_2) \in G\times G^{\rm op}$ is described below. Note that necessary, by the  definition of the partial action of $G\times G^{\rm op}$ on $K$, we have that $g_2=\gamma_1g_1\gamma_2$ for some $\gamma_1, \gamma_2 \in \Gamma.$
 
The algebra $C$ is generated by the linear span of elements of the form $M_{f_1}M_{f_2}M_\chi$ where $M_{f_1}$ is the multiplication operator with a continuous  function $f_1$ on $\pg$, $M_{f_2}$ is the multiplication operator with a continuous  function $f_2$ on  $\pg\o$ and $M_\chi$ where $\chi$ is the characteristic function of a coset of some subgroup $K_h$ for some $h\in G$. Since the action of  $G\times G^{\rm op}$ on $K$ is obviously as in the statement we will assume $\chi=1$.

Denote by $\lambda_g,\rho_g$ for $g$ in $G$ the left, respectively right convolution operators by $g$ on $\l^2(\Gamma)$. We have to compute
$$\lambda_{g_1}\rho_{g_2} M_{f_1}M_{f_2}(\lambda_{g_1}\rho_{g_2})^\ast=$$
$$[\lambda_{g_1}\rho_{g_2} M_{f_1}(\lambda_{g_1}\rho_{g_2})^\ast ] [\lambda_{g_1}\rho_{g_2}M_{f_2}(\lambda_{g_1}\rho_{g_2})^\ast].$$ 

In the Calkin algebra we have
$$\lambda_{g_1}\rho_{g_2} M_{f_1}(\lambda_{g_1}\rho_{g_2})^\ast=$$
$$\lambda_{g_1}\rho_{\gamma_1}\rho_{g_1}\rho_{\gamma_2} M_{f_1}(\lambda_{g_1}\rho_{\gamma_1}\rho_{g_1}\rho_{\gamma_2})^\ast=$$
$$\rho_{\gamma_1}\lambda_{g_1}\rho_{g_1}\rho_{\gamma_2}M_{f_1}(\rho_{\gamma_2})^\ast(\lambda_{g_1}\rho_{g_1})^\ast (\rho_{\gamma_1})^\ast=$$
$$\rho_{\gamma_1}\lambda_{g_1}\rho_{g_1}M_{f_1}(\lambda_{g_1}\rho_{g_1})^\ast (\rho_{\gamma_1})^\ast=$$
$$\rho_{\gamma_1}M_{f_1^{g_1}}(\rho_{\gamma_1})^\ast=M_{f_1^{g_1}},$$
where $f_1^{g_1}$ is the result of the action of $g_1$ on the function $f_1$.

 Note that in the computations above the  twisting cocycles do not intervene since the result of the conjugation by the unitaries implementing the crossed product, on the commutative algebra, in a twisted crossed product, does not depend on  the cocycle.
  
  It follows that the canonical action 
  of $G\times G\o$ on $\lc$ induces on $C$ the action described in the statement.     \end {proof}
  
  Using the natural action of $G$ on $\pg$ that we constructed above we deduce from Corollary \ref{p1} the following
  \begin{thm}\label{t3}
 The crossed product algebra 
 $$\mathcal D=C^{\ast}(( G \times  G^{\rm
op}) \ltimes_{\epsilon, \epsilon} C_0(\pg\times K \times \pg\o ))$$  coincides with the reduced
crossed product algebra
$$\mathcal D_{\rm red}= C^{\ast}_{\rm red}((G \times   G^{\rm
op}) \ltimes_{\epsilon, \epsilon} C_0(\pg\times K \times \pg\o )).$$
\end{thm}
  \begin{proof} We use the fact that the projection $\pi_d:\pg\rightarrow P^1(\mathbb R)$ (\cite{spiel}) is bijective except of the rational points where it is 2 to 1.
  
   We use Takesaki's Disintegration Theory (\cite {tak}, Chapter X, Theorem 3.8) for the crossed product algebra. Hence in any representation on a Hilbert space $H$ of the algebra $\mathcal D$, there exist a measure $\mu_{\mathcal D}$ on $\pg\times K \times \pg\o $ such that the Hilbert space $H$ becomes a field of Hilbert spaces over $\pg\times K \times \pg\o $ endowed with the measure $\mu_{\mathcal D}$. If $\mu_{\mathcal D}$
  gives mass 0 to the rational points then the result is exactly the content of Theorem \ref{t2}. Otherwise  the result still is a consequence  of Theorem \ref {t2}, which is this time used   for a double copy of the algebra in Theorem \ref{t2}, in  the case when the measure is supported on the rational points. Every such  measure is the sum of two singular measures for which Theorem \ref{t2} applies  and hence  the crossed product splits as a direct sum.
  
  \end{proof}
  
  We can now conclude the proof of Theorem \ref{t1}
  \begin{proof} 
 Indeed because of  Theorem \ref{extension} we can extend the representation $\Pi_Q$ of the algebra $\mathcal A$ to a representation $\tilde{\Pi_Q}$ of the algebra $\mathcal D=C^{\ast}(( G \times  G^{\rm
op}) \ltimes_{\epsilon, \epsilon} C_0(\pg\times K \times \pg\o ))$  into $Q(l^2(\Gamma))$. Because of Theorem \ref{t3} this factorises to the reduced
crossed product and hence also the restriction of $\tilde{\Pi_Q}$ to $\mathcal A$ factorises to the reduced crossed product algebra $\mathcal A\red$.
  
  \end{proof}

\section{Application to Hecke operators}

Recall from \cite{BC} that the algebra $\H_0$  of linearly generated by the  double cosets $[\Gamma\sigma\Gamma]$, for $\sigma \in G$ is represented into $B(l^2(\Gamma\backslash G))$, and correspondingly the reduced \csl- Hecke algebra $\H\red$ is the \csl-algebra generated by the image of $\H_0$. The canonical trace on $\H\red\subseteq B(l^2(\Gamma\backslash G))$  is implemented by the vector $[\Gamma]\in l^2(\Gamma\backslash G)$.

In \cite{ON} (see also \cite{Ra}) we proved that if $\pi_{13}$ is the 13-th projective unitary representation in the analytic series of representations of $\gr$ and $\zeta$ is a cyclic trace vector (see \cite {HGJ}, \cite{Jo})
then the following sum converges 
 $$t([\Gamma\sigma\Gamma])=\sum_{\theta \in \Gamma\sigma\Gamma}
 <\pi_{13}(\theta)\zeta, \zeta>\theta$$
 and defines an element in $\cs\red(G,\e)$.
 Moreover it is proven in \cite{ON} (see also \cite{Ra}) that
\begin{thm} \label{Ra}
(1)The application  mapping $[\Gamma\sigma\Gamma]$ into $t([\Gamma\sigma\Gamma])$, for $\sigma\in G$, extends by linearity to a trace preserving isomorphism from $\H\red$ into $\cs\red(G,\e)$.

(2)The   map  defined by linearly  extending   the formula

$$\Psi([\Gamma\sigma\Gamma])= E^{\mathcal L(G,\e)}_{\mathcal L(\Gamma,\e)}(t([\Gamma\sigma\Gamma])\cdot t([\Gamma\sigma\Gamma])), \sigma\in G$$
is a $\ast$-homeomorphism from $\H_0$ into  completely bounded maps on  $\mathcal L(G,\e)$. Here $E^{\mathcal L(G,\e)}_{\mathcal L(\Gamma,\e)}$ is the canonical type $II_1$  conditional expectation.

(3)If we extend $\Psi([\Gamma\sigma\Gamma])$ to $l^2(\Gamma)$ which is identified to the $L^2$-space associated to the II$_1$ factor $\mathcal L(\Gamma,\e)$, then $\Psi([\Gamma\sigma\Gamma])$ viewed as an element 
$\widetilde{\Psi([\Gamma\sigma\Gamma])}$ in  $B(l^2(\Gamma))$ is unitarily equivalent 
to the classical Hecke operator (\cite{He}) associated to $[\Gamma\sigma\Gamma]$ acting on Maass forms on the upper halfplane (modulo a commuting phase that is a function of the invariant laplacian).

(4) The linear application mapping $[\Gamma\sigma\Gamma]$ into
$$S([\Gamma\sigma\Gamma])=\chi_K(t([\Gamma\sigma\Gamma])\otimes 
t([\Gamma\sigma\Gamma]))\chi_K\in C^{\ast}((G \times G^{\rm
op}) \ltimes_{\epsilon, \epsilon} C(K)), $$
extends to a trace preserving isomorphism from $\H\red$ into $$\mathcal A\red=C^{\ast}\red((G \times G^{\rm
op}) \ltimes_{\epsilon, \epsilon} C(K)),$$  if the later algebra is endowed with the canonical trace, associated to the Haar measure on $K$. 
\end{thm}

As a corollary of the above theorem from \cite{ON} (see also \cite{Ra}) and 
because of Theorem \ref{t1} we obtain 

\begin{cor}\label{thecor} The linear application obtained by linear extension
by mapping $[\Gamma\sigma\Gamma]$  into the projection
 $\pi_Q( \widetilde{\Psi([\Gamma\sigma\Gamma])})$ into the Calkin algebra
 $Q(l^2(\Gamma))$, extends to a \csl-isomorphism from $\H\red$ with values into 
 $Q(l^2(\Gamma))$

\end{cor}

\begin{proof} 

Recall that $\pi_Q$ is the projection from $B(l^2(\Gamma))$ into the Calkin algebra, while $\Pi_Q$ is the representation of the algebra $\mathcal A$ obtained by composing the representation of $\mathcal A$ into $B(l^2(\Gamma))$ with the projection into the Calkin algebra.

Because of Theorem  \ref{t1} this will certainly hold true if we 
prove that in the Calkin algebra we have the equality:
\begin{equation}\label{yy}
\Pi_Q(S([\Gamma\sigma\Gamma]))= \pi_Q( \widetilde{\Psi([\Gamma\sigma\Gamma])}), \sigma \in G.
\end{equation}
Indeed, if the above  equality  holds true, then because of (4) in the previous theorem, the map $[\Gamma\sigma\Gamma]\rightarrow \pi_Q( \widetilde{\Psi([\Gamma\sigma\Gamma])}), \sigma \in G$ will extend to a continuous * homeomorphism from $\H\red$ into $Q(l^2(\Gamma))$. 

It remains to prove  the equality (\ref{yy}). By   definition, the selfadjoint operator $\widetilde{\Psi([\Gamma\sigma\Gamma])}$, in point (3) of Theorem \ref {Ra}, is obtained by extending $\Psi([\Gamma\sigma\Gamma])$ as a completely positive map on the von Neumann II$_1$ factor to an operator acting on $l^2(\Gamma)$. This  amounts to transform the conditional expectation
$E^{\mathcal L(G,\e)}_{\mathcal L(\Gamma,\e)}$ in point (2) of the above theorem into left and right multiplication by the characteristic function of $\Gamma$. This proves that the image of $S([\Gamma\sigma\Gamma])\in \mathcal A$ in the representation of the algebra $\mathcal A$ into $B(l^2(\Gamma))$ coincides with 
$\widetilde{\Psi([\Gamma\sigma\Gamma])}$ for $\sigma\in G$. 

\end{proof}

 Straightforwardly, from the previous corollary, we obtain that
\begin{cor}\label{Calkin}  For $\sigma \in G$,  the norm, in the Calkin algebra, of the Hecke operator  $\widetilde{\Psi([\Gamma\sigma\Gamma])}$  is equal to to the norm
of the Hecke operator associated to $[\Gamma\sigma\Gamma]$ in the reduced 
Hecke algebra $\H\red$ (i.e. the norm of $[\Gamma\sigma\Gamma]$ viewed as an operator acting on $l^2(\Gamma\backslash G)$).
\end{cor}

Let $T_{p^n}$ be the Hecke operator, acting on Maass forms, associated to 
$\sigma_{p^n}= \left(\begin{array}{cc} p^n & 0 \\ 0 & 1 \end{array} \right)$, $n$ a positive integer. As proved in part (3) in Theorem \ref{Ra}, this is unitarily equivalent  (modulo a phase operator that is a function of the invariant laplacian  on the upper half plane) to $\widetilde{\Psi([\Gamma\sigma_{p^n}\Gamma])}$.
The proof of Theorem \ref{thm3} is concluded as follows:

\ \begin{proof}  (Theorem \ref{thm3}). By part (3) of Theorem \ref{Ra} the essential  spectrum of $T_p$ coincides with the essential spectrum of $\widetilde{\Psi([\Gamma\sigma_{p}\Gamma])}$. Note the phase that intervenes in the unitarily equivalence of the two selfadjoint operators,, being a function of the laplacian  doesn't perturb the spectrum since the invariant laplacian commutes with all the Hecke operators.

By Corollary \ref {Calkin}, the spectrum, in the Calkin algebra of $\widetilde{\Psi([\Gamma\sigma_{p}\Gamma])}$ coincides with the spectrum of the  double coset $[\Gamma\sigma_{p}\Gamma]$ viewed as an operator acting on acting on $l^2(\Gamma\backslash G)$. But this is equal to the norm of the radial element $\chi_1$ (sum of words of length 1) in the free group $F_{(p+1)/2}$.
This is  because the $\H\red$ is identified (\cite{Kr}) with the radial algebra (\cite{TP}) in the free group. By \cite{Cohen} the spectrum of $\chi_1$, viewed as, an element of  $\cs\red(F_{(p+1)/2})$ is the interval $[-2\sqrt p, 2\sqrt p]$.
 \end{proof}

\section{Apendix}

In this appendix  we prove thew existence of  a "good" trace vector for the representation 
$\pi_{13}$ such  that the representation $t$ takes its values in $\cs(G,\e)$ (it is obvious that  it takes values in the von Neumann algebra associated to $G$ but it less obvious that it takes values in the \csl-algebra associated to $G$). This fact is not contained in \cite{ON}, but it is contained in the unpublished preprint \cite{Ra}. We use the notations introduced in the previous section.

\begin{lemma}
There exists a unit vector vector $\xi$ in $H_{13}$, that is a trace vector for the von Neumann algebra generated by the image of $\Gamma$ through the representation $\pi_{13}$ in $B(H_{13})$ and such that the representation $t$ of $\H\red$ associated to this representation takes values in $\cs(G,\e)$. Note that because   the von   Muray von Neumann dimension of $H_{13}$ (with respect to the above von Neumann algebra is 1) as proved in (\cite{HGJ}   this vector is also a cyclic vector.
\end{lemma}

{\it Proof.} Consider the space $\H_{13}$ of positive functions on $\PSL_2(\R)$ that are obtained
as matrix coefficients from elements $\eta$ in $H_{13}$ (that is $\varphi : G \to \C$ belongs
to $\H_{13}$ if there exists $\eta$ in $H_{13}$ such that $\varphi(g) = \langle \pi_{13}(g) \eta,\eta \rangle$,
$g$ in $\PSL_2(\R)$.

Obviously, $\H_{13}$ is a cone closed to infinite convergent sums. Indeed if $(\eta_i)$
is a family of vectors in $H_{13}$, $\sum\|\eta_i\|^2<\infty$, each determining the positive
functional $\varphi_i$. Consider the Hilbert subspace $L$ of $H_{13}\otimes \ell^2(I)$
generated by $\bigoplus\limits_{i \in I} \pi(g)\eta_i$. This space is
obviously invariant to the action of $G$. Since $\pi_{13}$ is irreducible $\pi(g)|_L$ is a multiple of 
the representation $\pi_{13}$ and because we have a cyclic vector, it is unitary equivalent to 
$\pi_{13}$. The vector $\eta = \bigoplus\limits_{i \in I} \eta_i$ will then determine the positive
definite function on $G$ defined by the formula
 $$\sum\limits_i \varphi_i(g) =\sum\limits_i \langle \pi_{13}(g) \eta_i,\eta_i \rangle, g \in G.$$

In the sequel we denote $\pi_{13}$ simply by $\pi$.

If choose a trace vector $\xi$ then by definition  $t^{\Gamma\sigma\Gamma}$ is equal to 
$$\sum\limits_{g\in [\Gamma\sigma\Gamma]}\overline{ \langle \pi(g) \xi,\xi \rangle }g.$$
 
If $\varphi_\eta(g) = \overline{\langle \pi(g)\eta,\eta \rangle}$, $g\in\PSL_2(\R)$ is determined by the vector 
$\eta$, then for $a$ in $L^1(\L(\Gamma, \e),\tau)$ the vector $\pi(a)\eta$
 will determine a functional $\varphi_a\in \H_{13}$, that is easily computed as 
\begin{equation}\label{equ}
\varphi_{\pi(a)\eta} | _{G} = a^* \varphi a.
\end{equation}

We are looking to find a positive functional in $\H_{13}$ that has the property that 
$\varphi | _{G}$ belongs to the reduced $C^*$-algebra of $G$, and such that moreover
$\varphi$ is implemented by a trace vector 


To find such a $\varphi=\varphi_\xi$ is therefore sufficient to find a vector $\xi$ such that the corresponding
positive functional has the following properties:

1) the restriction of $\varphi_\xi$ to $\Gamma\sigma\Gamma$ determines an element in
$C_{\rm red}^*(G, \varepsilon)$;

2) $\varphi_\xi|_\Gamma$ is invertible in $C_{\rm red}^*(\Gamma)$.

Indeed if we found such a vector $\xi$ then because of formula (\ref{equ}) it follows that the vector $\xi_0=
\pi((\varphi_\xi|_\Gamma)^{-1})\xi$ is a trace vector.

Moreover, let $t_0$ be the representation of the Hecke algebra associated to the vector $\xi_0$.
Then $$t_0^{\Gamma\sigma\Gamma}= (\varphi_\xi |_\Gamma)^{-1/2} \varphi_\xi |_{\Gamma\sigma\Gamma}
(\varphi_\xi |_\Gamma)^{-1/2}, $$ and this will then be an element   in $C_{\rm red}^*(G, \varepsilon)$  if $\xi$ is chosen as above.

We now use a result by  in  [BH] (proof of Theorem A1) which says that given $x\geq 0$, $x\neq 0$ in
$C_{\rm red}^*(\Gamma, \varepsilon)$ there exists unitaries $\gamma_1,\ldots,\gamma_n$ in $\Gamma$ such
that $\sum \gamma_i x \gamma_i^{-1}$ is invertible.

Let $\xi$ be a vector in $H_{13}$, generating a positive definite function on G,  which  has rapidly decreasing coefficients  so that property (1) holds. Note that to verify this it is sufficient to verify on every coset of $\Gamma$ in $G$ and this follows from the results in   (\cite {KS}). For example we may take the vector of evaluation
at $0$ in the model of the unit disk).

Then we construct the functional $\varphi_\xi$ and use the above mentioned result in [BH], to replace 
$\varphi_\xi$ by $\sum \gamma_i^{-1} \varphi_\xi \gamma_i = \Psi _0$.

Then $\Psi_0$ corresponds to the vector $\frac1{\sqrt n} ( \oplus\, \pi(\gamma_i)\xi)$
which is a vector generating a positive definite function with  rapidly decreasing coefficients. Consequently, by construction,  $\Psi_0|_\Gamma$ is invertible and hence its   inverse belongs to the $C^\ast$-algebra.


\begin{thebibliography}{LMD}
  
  

\bibitem{AO}  Akemann, C.A.,  Ostrand P.A.: On a tensor product $C^*$-algebra associated with the free group on two generators, J. Math. Soc. Japan. 27, 589-599, 1975.

\bibitem{CAD}  Anantharaman-Delaroche, C.:  Amenability and exactness for dynamical systems and their C$^\ast$ -algebras,  Trans. A.M.S.  354,  4153-4178, 2002.

\bibitem{BR} J. Bassi.,  R\u adulescu F.,  A mixing property for the action of $\SL(3,\mathbb{Z}) \times \SL(3,\mathbb{Z})$ on the Stone-Cech boundary of $\SL(3,\mathbb{Z})$, arXiv:2111.13885


\bibitem{BeCr}Alex Bearden and Jason Crann, Amenable dynamical systems over locally compact groups, Ergodic Theory Dynam. Systems (2021), DOI doi:10.1017/etds.2021.57

\bibitem{BH} Bridson, Martin R.; de la Harpe, Pierre Mapping class groups and outer automorphism groups of free groups are $C\sp *$-simple. J. Funct. Anal. 212 (2004), no. 1, 195--205. 





\bibitem{BC} J. B. Bost and A. Connes, Hecke algebras, type III factors and phase transitions with spontaneous symmetry breaking in number theory. Selecta Math. (N.S.) Vol.1 (1995), no. 3, 411-457.





\bibitem{BK}Jacek Brodzki, Chris Cave, and Kang Li, Exactness of locally compact groups, Adv. Math. 312 (2017), 209–233.

\bibitem{BO} Brown, N.,  Ozawa, N.:  C*-algebras and finite-dimensional approximations.
Graduate Studies in Mathematics, 88. American Mathematical Society, 2008.





\bibitem{BE} A. Buss, S. Echterhoff, R. Willett, Amenability and weak containment for actions of locally compact groups on \csl-algebras, arXiv:2003.03469




\bibitem
{Ca} Calkin, J.W.: Two-sided ideals and congruences in the ring of bounded operators in Hilbert space. 
Ann. of Math. 42,  839-873, 1941.

\bibitem{Cohen}J.M. Cohen: Operator norms on free groups, Boll. Un. Mat. Ital. B (6) 1 (1982), no. 3, 1055-1065



\bibitem{cooper} Cooper, Automorphisms of free groups have finitely generated fixed point sets, Journal of Algebra
Volume 111, Issue 2, December 1987, Pages 453-456








 

\bibitem{Cut} Cutland, N.:  Loeb measures in practice: recent advances. Lecture Notes in Mathematics, 1751. Springer-Verlag, Berlin, 2000.


\bibitem{GH} de la Harpe, Pierre; Ghys, Etienne (1990), Sur les groupes hyperboliques d'après Mikhael Gromov (in French), Birkhäuser







\bibitem{TP} A. Figa-Talamanca, M.A. Picardello, Harmonic analysis on free groups, Lecture Notes in Pure and Appl. Math. 87, Dekker, New York, 1983.

\bibitem{F1} H. Furstenberg, Random walks and discrete subgroups of Lie groups, 1971. Advances
in Probability and Related Topics, Vol. 1, pp. 1–63, Kekker, New York.

\bibitem{F2} H. Furstenberg. Boundary theory and stochastic processes on homogeneous spaces. Harmonic analysis on homogeneous spaces (Proc. Sympos. Pure Math., Vol. XXVI, Williams Coll., Williamstown, Mass., 1972). 1973.

\bibitem{GG} I.M. Gelfand, M.I. Graev, I. Piatetsky-Shapiro
Representation Theory and Automorphic Functions
Saunders, San Francisco (1969)

\bibitem{HGJ} F.M. Goodman, P. De La Harpe, V.F.R. Jones, Coxeter graphs and towers of algebras, Springer-Verlag, New York, 1989. 


\bibitem{Gr} R. Grigorchuk, Thesis

\bibitem{gkn} R. I. Grigorchuk, V. A. Kaimanovich, and T. Nagnibeda, Ergodic properties of boundary actions and Nielsen-Schreier theory, Advances in Mathematics 230
(3) (2012), p. 1340-1380.




\bibitem{GHW} Gunter, E., Higson, N.:, Weinberger, S.,  The Novikov Conjecture for Linear Groups, Publications mathématiques de l'IHÉS volume 101, pages 243–268 (2005)








\bibitem{GH}  Guentner, E.,  Higson, N.: Group 
$C^*$-algebras and $K$-theory. Noncommutative geometry, 137--251, Lecture Notes in Math., 
1831, Springer, Berlin, 2004.


\bibitem{Ha}Haagerup, U., The standard from of von Neumann algebras. Math. Scand. 37 (1975), 271-283.

\bibitem{He} Hecke, E.:  Uber Modulfunktionen und die Dirichletschen Reihen mit Eulerscher Produktentwicklung. I., Mathematische Annalen, 114: 1-28, (1937).




\bibitem{Ik}Akio Ikunishi, The W$^\ast$-dynamical system associated with a C$^\ast$-dynamical system, and un-
bounded derivations, J. Funct. Anal. 79 (1988), no. 1, 1–8

\bibitem{Jo} V. Jones, Bergman space zero sets, modular forms, von Neumann algebras and ordered groups, arXiv:2006.16419



\bibitem{kks} Vadim Kaimanovich, Ilya Kapovich, Paul Schupp, The Subadditive Ergodic Theorem and generic stretching factors for free group automorphisms, Israel J. Math. 157 (2007), pp. 1-46 

\bibitem{Klep}
A.\ Kleppner,
\emph{The structure of some induced representations}.
Duke Math.\ J.\ \textbf{29} (1962), 555--572.









\bibitem{Kr} A. Krieg; Hecke algebras. Mem. Amer. Math. Soc. 87 (1990), no. 435, x+158 pp.

\bibitem{KS} Kuhn, Gabriella; Steger, Tim. More irreducible boundary representations of free groups. Duke Math. J. 82 (1996), no. 2, 381--436.


\bibitem{Li} Lindstrom, T.: An invitation to nonstandard analysis, in Nonstandard analysis and its applications, Hull, 1986, 1-105, London Math. Soc. Stud. Texts, 10, Cambridge Univ. Press, Cambridge, 1988.

\bibitem{Lo} Loeb, P.A.: Conversion from nonstandard to standard measure spaces
and applications in probability theory, Trans. Amer. Math. Soc.  211,
113-122, 1975.



 





\bibitem{Mats} Masayoshi Matsumura, A characterization of amenability of group actions on C$\ast$-algebras,
J. Operator Theory 72 (2014), no. 1, 41–47, DOI 10.7900/jot.2012sep07.1958


\bibitem{Oz1} Ozawa, N.:  Amenable actions and exactness for discrete groups. 
C. R. Acad. Sci. Paris S\'er. I Math. 330,  691Ð695,  2000.


\bibitem{Oz} Ozawa, N.: Solid von Neumann algebras. Acta Math. 192, 111-117, 2004.

\bibitem{Oz2} Ozawa, N.: Amenable actions and applications, Proceedings I. C. M., Madrid, Spain, 2006,
E. M. S., 1563-1580, 2006.

\bibitem{OS}Narutaka Ozawa and Yuhei Suzuki, On characterizations of amenable C$^\ast$-dynamical systems and new examples, Selecta Math. (N.S.) 27 (2021), no. 5, Paper No. 92, 29

\bibitem{pac} Judith Packer, I Raeburn, Twisted crossed products of C*-algebras, Mathematical Proceedings of the Cambridge Philosophical Society, Cited by 81
Get access  Volume 106, Issue 2 September 1989 , pp. 293-311

\bibitem{Pe} Pedersen, G. C$^\ast$-algebras and their automorphism groups, London Math.  Soc. Monographs 14, 
Academic Press, 1979.




\bibitem{ON}F. R\u adulescu, Free group factors and Hecke operators, notes taken by N. Ozawa, Proceedings of the 24th Conference in Operator Theory, Theta Advanced Series in Mathematics, Theta Foundation, 2014.


\bibitem{Ra} F. R\u adulescu, Type II$_1$ von Neumann representations for Hecke operators on Maass forms 
and Ramanujan-Petersson conjecture, preprint arxiv arXiv:0802.3548, 2008.

\bibitem{Ra1} F. R\u adulescu, Endomorphisms of spaces of virtual vectors fixed by a discrete group, Russ. Math. Surv. 71 291-343, 2016

\bibitem{Ra2} F. R\u adulescu, The operator algebra content of the Ramanujan–Petersson problem, J. Non-commutative Geometry, 13, 805-855, 2019.



\bibitem{raney}Raney, G. N. On continued fractions and finite automata. Math. Annalen 206, (1973) 265-283.







\bibitem{JR} J. Renault, A Groupoid Approach to $\cs$-Algebras,  Lecture Notes in Mathematics, volume 793, Springer Verlag, 1980.







\bibitem{Sa} Sako, H.:  The class 
$S$ as an $ME$ invariant. Int. Math. Res. Not.  15,   2749-2759, 2009.

\bibitem{Sal} Sally, An Introduction to p-adic Fields, Harmonic Analysis and the Representation Theory of SL2, Letters in Mathematical Physics 46: 1-47, 1998.

\bibitem{Serre} J.P.  Serre, A Course in Arithmetic, Springer Verlag, 1973/

\bibitem{series} C. Series,  Non-Euclidean geometry, continued fractions, and ergodic theory. Math. Intelligencer 4 (1982), 24-31. 

\bibitem{Sk} Skandalis, G.:
Une notion de nucl\'earit\'e en  $K$-th\'eorie (d'apr\`es J. Cuntz).  
[A notion of nuclearity in $K$-theory (following J. Cuntz)] K-Theory 
1, 549-573, 1988. 


\bibitem{spiel}J. Spielberg, Cuntz Krieger algebras associated with Fuchsian groups Ergodic Theory Dynamical Systems, Volume 13, Issue 3 September 1993 , pp. 581-595






Tunbridge Wells, 1981.



\bibitem{tak} Takesaki, M., Theory of operator algebras. II, Encyclopaedia of Mathematical Sciences volume 125

\bibitem{unt} Unterberger, A.,  The Ramanujan-Petersson conjecture for Maass forms, arXiv:2001.10956.

\bibitem{Ze}H. Zettl, A characterization of ternary rings of operators. Adv. Math. 48 (1983), 117-143.






\end{thebibliography}
\end{document}